\documentclass{amsart}
\usepackage{amsmath,amsthm,amssymb,times,graphics}

\theoremstyle{plain}

\newtheorem*{theorem*}{Theorem}

\newtheorem*{lemma*}{Lemma}

\theoremstyle{definition}

\newtheorem{problem}{Problem}

\newtheorem*{answer}{Answer}

\newtheorem*{definition*}{Definition}

\newtheorem*{proposition*}{Proposition}

\newcommand{\<}{\langle}

\renewcommand{\>}{\rangle}

\newcommand{\Z}{\mathbb{Z}}

\newcommand{\A}{\mathcal{A}}
\newcommand{\power}{2^{\{P_1, \cdots , P_n\}}}

\newcommand{\bp}{\begin{problem}}

\newcommand{\ep}{\end{problem}}

\newcommand{\ba}{\begin{answer}}

\newcommand{\ea}{\end{answer}}

\newcommand{\ben}{\renewcommand{\theenumi}{\alph{enumi}}

\renewcommand{\labelenumi}{(\theenumi)}\begin{enumerate}}

\newcommand{\een}{\end{enumerate}}

\begin{document}

\title{Secret Sharing using Non-Commutative Groups and the Shortlex Order}
\author{Bren Cavallo}
\address{CUNY Graduate Center, City University of New York}
\email{bcavallo@gc.cuny.edu}

\author{Delaram Kahrobaei}
\address{CUNY Graduate Center and City Tech, City University of New York}
\email{dkahrobaei@gc.cuny.edu}


\begin{abstract}
In this paper we review the Habeeb-Kahrobaei-Shpilrain secret sharing scheme \cite{DBLP:journals/corr/abs-1205-0157} and introduce a variation based on the shortlex order on a free group.  Drawing inspiration from adjustments to classical schemes, we also present a method that allows for the protocol to remain secure after multiple secrets are shared.
\end{abstract}

\maketitle

\section{Introduction}
\nocite{*}
Secret sharing is a cryptographic protocol by which a dealer distributes a secret via shares to participants such that only certain subsets of participants can together use their shares to recover the secret. A secret sharing scheme begins with a dealer, a secret, participants, and an access structure. The access structure determines which groups of participants have access to the secret. The goal of the scheme is to distribute the secret to the participants in such a way that only sets of participants within the access structure have access to the secret. In this way, it is most often the case that no individual participant can recover the secret on their own.

\

Secret sharing schemes are ideal tools for when the secret is both highly important and highly sensitive. The fact that there are multiple shares, as opposed to one private key in private key cryptography, makes the secret less likely to be lost while allowing high levels of confidentiality. If any one share is compromised the secret can generally still be recovered with the non-compromised shares. Additionally, even though the secret is spread out over multiple shares, recovering the secret is limited by the access structure, and so the secret remains secure. Secret sharing has applications in multi-party encryption, Byzantine agreement, and threshold encryption among others. See \cite{Beimel:2011:SSS:2017916.2017918} for a survey on secret sharing and its applications in cryptography and computer science.
\section{Formal Definition}
A $\emph{secret sharing scheme}$ consists of a dealer, $n$ participants, $P_1, \ldots P_n$, and an access structure $\A \subseteq \power$ such that for all $A \in \A$ and $A \subseteq B$, $B \in \A$.
\
To share a secret $s$, the dealer runs an algorithm:
$$Share(s) = (s_1, \cdots, s_n)$$ and then distributes each share $s_i$ to $P_i$.
\
In order to recover the secret, participants can run the algorithm $\emph{Recover}$ which has the property that for all $A \in \A$:
$$Recover(\{s_i : i \in A\}) = s$$
and if $A \notin \A$ then running $\emph{Recover}$ is either computationally infeasible or impossible.
\
As such, only groups of participants in $\A$ can access the secret. The monotonicity of $\A$ is also apparent in that if $A \in \A$ and $A \subseteq B$ then the set of participants in $A$ could also recover the secret for $B$. A secret sharing scheme is called $\emph{perfect}$ if $\forall A \notin \A$ the shares $s_i \in A$ together give no information about $s$.
\section{Shamir's Secret Sharing Scheme}
One of the more common access structures one sees in secret sharing is the $\emph{(k,n) threshold}$: $$\A = \{A \in \power : |A| \geq k\}$$
Namely, $\A$ consists of all subsets of the $n$ participants of size $k$ or greater. We call a secret sharing scheme that has $\A$ as a $(k,n)$ threshold a $\emph{(k,n) threshold scheme}$. The problem of discovering a perfect $(k,n)$ threshold scheme was solved independently by G. Blakely \cite{Blakley1979} and A. Shamir \cite{DBLP:journals/cacm/Shamir79} in 1979.

\

In the Shamir Secret Sharing Scheme, the secret is an element in $\Z_p$ where $p$ is a prime number larger than the number of participants. Given a secret $s$, the dealer generates the shares for a $(k,n)$ threshold by doing the following:
\begin{itemize}
\item The dealer randomly selects $a_1, \cdots, a_{k-1}\in \Z_p$ such that $a_{k-1} \neq 0$ and constructs the polynomial $f(x) = a_{k-1}x^{k-1}+\cdots+a_1 x + s$
\item For each participant $P_i$ the dealer publishes a corresponding $x_i \in \Z_p$. The dealer then distributes the share $s_i = f(x_i)$ to each $P_i$ over a private channel.
\end{itemize}
Any subset of $k$ participants can then reconstruct the polynomial $f(x)$ by using polynomial interpolation and then finding $f(0) = s$. This method finds $s$ uniquely as any degree $k-1$ polynomial is uniquely determined by the $k$ shares. The shares are consistent because each $(x_i, f(x_i))$ is a point on the polynomial $f(x)$ and thus any $k$ shares will reconstruct the same polynomial. In order to reconstruct a polynomial $f(x) = a_0 + a_1 x + \cdots a_{k-1} x^ {k-1}$ given points $(x_1, f(x_1)), \cdots , (x_k, f(x_k))$ one can solve for the coefficients column in the following system of linear equations:
$$\begin{pmatrix}
{x_1}^{k-1} & \cdots & x_1 & 1 \\
{x_2}^{k-1} & \cdots & x_2 & 1 \\
\vdots & \ddots & \vdots & \vdots \\
{x_k}^{k-1} & \cdots & x_k & 1 \end{pmatrix}
\begin{pmatrix}
a_{k-1} \\ a_{k-2} \\ \vdots \\ a_0
\end{pmatrix} =
\begin{pmatrix}
f(x_1) \\ f(x_2) \\ \vdots \\ f(x_k)\end{pmatrix}$$
\
The above method of interpolation demonstrates that Shamir's scheme is perfect. If there were less than $k$ shares, than the system of equations above would have more equations than unknowns, and there would not be a unique solution for $a_0$.
\section{Secret Sharing Using Non-commutative Groups}
Given a set of letters $X = \{x_1, x_2, \ldots, x_n\}$ we define the free group generated by $X$, $F(X)$, as the set of reduced words in the alphabet $X^{\pm 1} = \{x_1^{\pm 1}, \ldots , x_n^{\pm 1}\}$, where a word is reduced if there are no subwords of the form $x_i^{-1} x_i$ or $x_i x_i^{-1}$. Given a set of words $R \subset F(X)$ we define $\<\<R\>\>$ as the smallest normal subgroup of $F(X)$ containing $R$ and define the group $G = \< X | R \> = F(X)/\<\<R\>\>$. We call $R$ the set of $\emph{relators}$ of $G$.

\

A group $G = \<X | R\>$ has a \emph{solvable word problem} if there exists an algorithm to determine if any word $w \in G$ is trivial. Habeeb-Kahrobaei-Shpilrain (HKS) secret sharing \cite{DBLP:journals/corr/abs-1205-0157} uses a group with an efficiently solvable word problem to create an $(n,n)$ threshold scheme which can be extended to a $(k,n)$ threshold scheme using the method of Shamir.
\subsection{$(n,n)$ Threshold}
In this case the secret, $s$, is an element of $\{0,1\}^k$ which we view as a column vector. The setting is initialized by making a set of generators $X = \{x_1, \cdots , x_n\}$ public. To distribute the shares the dealer does the following:
\begin{itemize}
\item Distributes to each $P_i$ over a private channel a set of words $R_i$ in the alphabet $X^{\pm 1}$ that define the group $G_i = \<X | R_i\>$.
\item Randomly generates the shares $s_i \in \{0,1\}^k$ for $i=1,\cdots,n-1$ and $s_n = s - \sum_{j=0}^{n-1} s_j$ where the addition is bitwise addition in $\mathbb{F}_2^k$.
\item Publishes words $w_{ji}$ over the alphabet $X^{\pm 1}$ such that a word $w_{ji}$ is trivial in $G_i$ if $s_{ji} = 1$ and non-trivial if $s_{ji} = 0$.
\end{itemize}
Since the $G_i$ have efficiently solvable word problem, the participant $P_k$ can determine which of the $w_{jk}$ are trivial or non-trivial and can independently recover $s_k$. To recover the secret, the $P_i$ add the $s_i$ and find $s$. Note that even though the $w_{ji}$ are sent over an open channel, the shares remain secure since the $R_i$ are private. Therefore no other participant can recover $s_i$ from the $w_{ji}$ since only $P_i$ knows $G_i$.
\subsection{$(k,n)$ Threshold}
One can extend the above scheme to a $(k,n)$ threshold via Shamir's scheme. As is the case with Shamir's scheme, the secret $s$ is an element of $\Z_p$ and the shares, $s_i$, correspond to points on a polynomial of degree $k-1$ with constant term $s$. The shares are distributed and reconstructed in an identical manner as above by viewing the $s_i$ in their binary form. The trivial and non-trivial words are sent to each $P_i$ so that they reconstruct each $s_i$ in its binary form. After recovering their shares any element of the access structure can use polynomial interpolation to find $s$:
\begin{itemize}
\item The dealer randomly selects $a_1, \cdots, a_{k-1}\in \Z_p$ such that $a_{k-1} \neq 0$ and constructs the polynomial $f(x) = a_{k-1}x^{k-1}+\cdots+a_1 x + s$.
\item For each participant $P_i$ the dealer publishes a corresponding $x_i \in \Z_p$. The dealer then converts each $s_i = f(x_i)$ into binary. And thus, each $s_i$ can be viewed as a column vector of length $l = \left \lfloor {\log_2 {p}} \right \rfloor + 1$.
\item As was the case in the $(n,n)$ scheme, the dealer distributes the $s_i$ over an open channel by sending each $P_i$ the words $w_{1i}, \cdots, w_{li}$ over the alphabet $X^{\pm}$ such that $w_{ji}$ is trivial in $G_i$ if $s_{ji} = 1$ and non-trivial if $s_{ji} = 0$.
\item The participants reconstruct their own $s_i$ and can recover the secret using polynomial interpolation.
\end{itemize}
Some advantages this secret sharing scheme has over Shamir's scheme include the fact that after the $R_i$ are distributed, one can still use them to send out and reconstruct more secrets rather than having to privately distribute new shares each time a different secret is picked. Private information has to only be sent once initially for an arbitrary amount of secrets to be shared due to the method of distributing the shares. Despite this, the scheme is vulnerable to an adversary determining the relators by seeing patterns in words they learn are trivial. Namely, after a participant reveals their share (possibly while recovering the secret) an adversary could determine which of the $w_{ji}$ were trivial and potentially find the group presentation of $G_i$ which would allow them to reconstruct $P_i$'s share on their own. As in \cite{DBLP:journals/corr/abs-1205-0157}, we assume that this is a computationally difficult problem. Moreover, in {\bf Section 5} we provide a method to update relators over time thus limiting the amount of information an adversary could obtain about a single group. Another advantage to this scheme is that since it is based on the Shamir secret sharing protocol it can benefit from the large amount of research done on Shamir's scheme. For instance, the verification methods or proactive secret sharing protocols from \cite{DBLP:conf/eurocrypt/Stadler96} and \cite{Herzberg:1997:PPK:266420.266442} can still be used in this scheme.
\subsection{Small Cancellation Groups}
In this section we introduce a candidate group for the above secret sharing scheme.

\

A word $w$ is \emph{cyclically reduced} if it is reduced in all of its cyclic permutations. Note that this only occurs if the word is freely reduced, it has no subwords of the form $x_i^{-1} x_i$ or $x_i x_i^{-1}$, and the first and last letters are not inverses of each other.

\

A set of words $R$ is called \emph{symmetrized} if each word is cyclically reduced and the entire set and their inverses are closed under cyclic permutation. If $R$ is viewed as a set of relators, symmetrizing $R$ does not change the resulting group as the closure $R$ under cyclic permutations and inverses is a subset of the normal closure.

\

Given a set $R$ we say that $v$ is a \emph{piece} if it is a maximal initial subword of two different words, namely if there exist $w_1, w_2 \in R$ such that $w_1 = vr_1$ and $w_2 = vr_2$. A group $G = \<X | R\>$ satisfies the \emph{small cancellation condition $C'(\lambda)$} for $0 < \lambda < 1$ if for all $r \in R$ such that $r = vw$ where $v$ is a piece, then $|v| < \lambda |r|$.

\

Small cancellation groups satisfying $C'(\frac{1}{6})$ have a linear time algorithm for the word problem \cite{domanski:the} making them an ideal candidate for the HKS secret sharing scheme. Moreover, it can be seen from their definition that if the number of generators is large compared to the number of relators and lengths of the relators, it is likely that there will be small cancellation since the probability that any two words have a large maximal initial segment is low. After generating a random set of relators satisfying the above properties, it is also fast to symmetrize the set and then find the pieces and check that they are no larger than one sixth of the word. As such, it is fast to create such groups by repeatedly randomly generating relators, symmetrizing, and checking to see if they satisfy the $C'(\frac{1}{6})$ condition. There are other groups that have an efficient word problem that could also function as candidate groups, but small cancellation groups have the advantage of being efficient to generate randomly.
\subsection{Secret Sharing and the Shortlex Ordering}
Let $X = \{x_1, \cdots ,x_n\}$ and $G = \< X \>$. A shortlex ordering on G is induced by an order on $X^{\pm 1}$ as follows. Given reduced $w = x_{i_1} \cdots x_{i_p}$ and $l = x_{j_1} \cdots x_{j_k}$ with $w \neq l$ then $w < l$ if and only if:
\begin{itemize}
\item $|w| < |l|$
\item or if $p = k$ and $x_{i_a} < x_{j_a}$ where $a = \min_\alpha \{x_{i_\alpha} \neq x_{j_\alpha}\}$
\end{itemize}
For example, let $X = \{x,y\}$ and give $X^\pm$ the ordering $x < x^{-1} < y < y^{-1}$. Then some of the first words in order would be:

\

$e < x < x^{-1} < y < y^{-1} < x^2 < xy < xy^{-1}
< x^{-2} < x^{-1} y < x^{-1}y^{-1}<yx<yx^{-1}< y^2 < y^{-1}x<y^{-1}x^{-1}<y^{-2}<x^3< x^2 y < x^2 y^{-1}
< xyx < xyx^{-1}<\cdots$

\

Utilizing the the shortlex ordering, we can modify the HKS $(k, n)$ threshold as follows:
\begin{itemize}
\item The dealer publishes the letters $X$ and over a private channel sends a set of words, $R_i$ in $X^{\pm 1}$ to each $P_i$ such that $G_i = \< X | R_i\>$ is a group with an efficient algorithm to reduce words with respect to the $R_i$  or compute normal forms.
\item The dealer chooses a secret $s \in \Z_{p}$ for some large prime $p > n$ and generates a random polynomial, $f$ in $\Z_p[x]$ with constant term $s$
\item The dealer assigns a public $x_i \in \Z_p$ to each participant, computes $f(x_i)$, and finds $s_i \in F(X)$ such that $s_i$ is the $f(x_i)^{th}$ word in $F(X)$. Note that $x_i$ is not a generator of $G$, but rather the $x$-coordinate associated to each participant's share.
\item The dealer publishes a word $w_i$ that reduces to $s_i$ in $G_i$. This can be done efficiently by interspersing conjugated products of relators between the letters of $s_i$.
\item Each participant $P_i$ computes their share by reducing $w_i$ to get $s_i$ and then computing its position in $F(X)$.
\item Using their shares they find the secret using polynomial interpolation.
\end{itemize}

\

The main advantage of this new method is that participants need only reduce one word rather than a number of words corresponding to the length of the secret. In general, being able to reduce words is more general than being able to solve the word problem in a finitely presented group and in some cases may be more complex.  It is important to note the following about this scheme:
\begin{itemize}
\item Given an algorithm that reduces words, each $w_i$ must reduce uniquely to $s_i$. This implies that if our reduction algorithm does not terminate at $s_i$, then it is not a viable share for this scheme. In that case, if a random $f(x_i)$ does not correspond to a fully reduced word or a word in normal form, the dealer can always assign $P_i$ a different $x_i$. It may also be necessary to check that each $w_i$ reduces to $s_i$ give the reduction algorithm before the shares are distributed.
\item Some reduction algorithms can be done in multiple ways given the same initial conditions and can terminate at different words.  As such, it is important to fix a protocol so that whatever process $P_i$ uses to reduce $w_i$ terminates at $s_i$.
\end{itemize}
\
\subsection{Platform Group}
For this variant of the HKS secret sharing scheme, we also propose $C'(\frac{1}{6})$ groups. Additionally, we propose the parameters $|X| = 40$, $|R| = 4$, and $|r| = 9$ for all $r \in R$. We find that with such parameters, generating a single $C'(\frac{1}{6})$ group can be done in roughly 1 second in GAP \cite{GAP4} by generating random relators of the given length and then checking that the set of relators satisfies the small cancellation condition. In order to reduce the $w_i$ to $s_i$, participants can use Dehn's algorithm which terminates in linear time \cite{domanski:the}. It is not guaranteed in general that Dehn's algorithm will reduce each $w_i$ to $s_i$, as such it is necessary to check that each $w_i$ reduces to $s_i$.
\
In order to test the efficacy of Dehn's algorithm in $C'(\frac{1}{6})$ groups for the purposes of this secret sharing scheme, we performed the following tests in GAP \cite{GAP4}:
\begin{itemize}
\item Generate 10 small cancellation groups using the parameters from the first paragraph of this section.
\item In each group we generated 100 words of length less than 10 and created corresponding large unreduced words of length $~500$ by inserting
conjugated products of relators between letters in our original word.
\item Applied an implementation of Dehn's algorithm due to Chris Staecker \cite{gapdehn} and checked that our unreduced word successfully reduced to the original word.
\end{itemize}
After running said tests, we found that Dehn's algorithm successfully reduced every word. The size considerations in the second item were given in part because there are enough non-trivial, Dehn reduced, words of length 10 or less in the free group on 40 generators to be used as shares in a practical setting.

\subsection{Efficiency}
Each step in modified HKS scheme can be done efficiently. As mentioned previously, generating $C'(\frac{1}{6})$ groups can be done quickly by repeatedly generating sets of relators and checking to see if they satisfy the necessary small cancellation condition. The necessary computations using the shortlex ordering can be done using basic combinatorial formulas that are very fast for a computer to evaluate. Additionally, the $w_i$ can be created efficiently from the $s_i$ by inserting conjugated products of relators and then reduced in polynomial time using Dehn's algorithm. Moreover, the dealer can also check that the $w_i$ reduce to the $s_i$ efficiently. Hence each additional step to the standard Shamir's scheme can be done efficiently. This is also an improvement over the standard HKS scheme since the amount of words that need to be reduced is independent of the length of the secret, making it possible for larger secrets to be distributed efficiently.
\section{Updating Relators}
The main security concern for this cryptoscheme is the possibility of an adversary discovering a participant's set of relators. This can either be done using information gained from combining shares, but even potentially just from the public $w_i$. As more secrets are shared, the original set of relators becomes less secure. Moreover, information may be discovered either by breaking into wherever a participant stores their relators or if partial information was discovered during the initial step. In this section we present a method to refresh a participant's relator set using the same inherent security assumptions necessary for the cryptoscheme, namely that at least one round of secret sharing is secure.  To do this we add steps that can take place before any new secret is sent out:
\begin{itemize}
\item For each $P_i$ the dealer creates a set of words, $R'_i$, in $X^{\pm 1}$ such that $G_i = \<X | R'_i \>$ satisfies the same desired properties.
\item In order to distribute each $r \in R'_i$, the dealer pads $r$ with relators in $R_i$ as done previously and publishes them.
\item $P_i$ then reduces $r$ by using the relators in $R_i$.
\item After the full set of words in $R'_i$ is published and reduced, $P_i$ deletes the original $R_i$ and sets $R_i := R'_i$.
\end{itemize}
If these steps are done before an adversary can gain adequate information about relators, then after an update phase the information an adversary has gained will be largely rendered useless. Also note that a single secret can be kept secure over a long period of time using the methods in \cite{Herzberg:1997:PPK:266420.266442}. In this case, it is important that the words in $R'_i$ are reduced with respect to the original $R_i$. As such, $R_i$ and $R'_i$ are not completely unrelated, but as the relators become updated each additional time, they will have less and less to do with the original set of relators.
\section{Conclusion}
In this paper we propose a modification of the HKS secret sharing scheme using the shortlex ordering on free groups. It improves the original scheme by removing the relation of the number of times each participant has to solve the word problem to the length of the secret. As such, larger secrets can be shared efficiently and the overall scheme is more efficient. Moreover, it shares the advantage over Shamir's scheme that multiple secrets can be shared given the same initial private information. We also introduce a method to update relators so that the scheme remains secure when arbitrarily many secrets are shared and that does not involve more private information being distributed.
\subsection{Support}
Delaram Kahrobaei is partially supported by the Office of Naval Research grant N00014120758 and also supported by PSC-CUNY grant from the CUNY research foundation, as well as the City Tech foundation.

\bibliographystyle{plain}
\bibliography{biblio}

\begin{thebibliography}{10}

\bibitem{Beimel:2011:SSS:2017916.2017918}
Amos Beimel.
\newblock Secret-sharing schemes: a survey.
\newblock In {\em Proceedings of the Third international conference on Coding
  and cryptology}, IWCC'11, pages 11--46, Berlin, Heidelberg, 2011.
  Springer-Verlag.

\bibitem{Blakley1979}
G.R. Blakley.
\newblock Safeguarding cryptographic keys.
\newblock In {\em Proceedings of the 1979 AFIPS National Computer Conference},
  pages 313--317, Monval, NJ, USA. AFIPS Press.

\bibitem{domanski:the}
B.~Domanski and M.~Anshel.
\newblock The complexity of dehn's algorithm for word problems in groups.
\newblock {\em J. Algorithms}, pages 543--549, 1985.

\bibitem{DBLP:conf/crypto/1991}
Joan Feigenbaum, editor.
\newblock {\em Advances in Cryptology - CRYPTO '91, 11th Annual International
  Cryptology Conference, Santa Barbara, California, USA, August 11-15, 1991,
  Proceedings}, volume 576 of {\em Lecture Notes in Computer Science}.
  Springer, 1992.

\bibitem{Feldman:1987:PSN:1382440.1383000}
Paul Feldman.
\newblock A practical scheme for non-interactive verifiable secret sharing.
\newblock In {\em Proceedings of the 28th Annual Symposium on Foundations of
  Computer Science}, SFCS '87, pages 427--438, Washington, DC, USA, 1987. IEEE
  Computer Society.

\bibitem{GAP4}
The GAP~Group.
\newblock {\em {GAP -- Groups, Algorithms, and Programming, Version 4.7.6}},
  2014.
\newblock http://www.gap-system.org.

\bibitem{DBLP:journals/corr/abs-1205-0157}
Maggie Habeeb, Delaram Kahrobaei, and Vladimir Shpilrain.
\newblock A secret sharing scheme based on group presentations and the word
  problem.
\newblock {\em Contemp. Math., Amer. Math. Soc.}, 582:143-- 150, 2012.

\bibitem{Herzberg:1997:PPK:266420.266442}
Amir Herzberg, Markus Jakobsson, Stanisllaw Jarecki, Hugo Krawczyk, and Moti
  Yung.
\newblock Proactive public key and signature systems.
\newblock In {\em Proceedings of the 4th ACM conference on Computer and
  communications security}, CCS '97, pages 100--110, New York, NY, USA, 1997.
  ACM.

\bibitem{handbook}
Derek~F. Holt, Bettina Eick, and Eamonn~A. O'Brien.
\newblock {\em Handbook of Computational Group Theory}.
\newblock CRC Press, 2005.

\bibitem{jarecki1996proactive}
S.M. Jarecki.
\newblock {\em Proactive Secret Sharing and Public Key Cryptosystems}.
\newblock Massachusetts Institute of Technology, Department of Electrical
  Engineering and Computer Science, 1996.

\bibitem{DBLP:books/crc/KatzLindell2007}
Jonathan Katz and Yehuda Lindell.
\newblock {\em Introduction to Modern Cryptography}.
\newblock Chapman and Hall/CRC Press, 2007.

\bibitem{DBLP:conf/eurocrypt/96}
Ueli~M. Maurer, editor.
\newblock {\em Advances in Cryptology - EUROCRYPT '96, International Conference
  on the Theory and Application of Cryptographic Techniques, Saragossa, Spain,
  May 12-16, 1996, Proceeding}, volume 1070 of {\em Lecture Notes in Computer
  Science}. Springer, 1996.

\bibitem{groupcrypto}
Alexei Myasnikov, Vladimir Shpilrain, and Alexander Ushakov.
\newblock {\em Group-based Cryptography}.
\newblock Springer, 2008.

\bibitem{DBLP:conf/crypto/Pedersen91}
Torben~P. Pedersen.
\newblock Non-interactive and information-theoretic secure verifiable secret
  sharing.
\newblock In Feigenbaum \cite{DBLP:conf/crypto/1991}, pages 129--140.

\bibitem{DBLP:journals/cacm/Shamir79}
Adi Shamir.
\newblock How to share a secret.
\newblock {\em Commun. ACM}, 22(11):612--613, 1979.

\bibitem{DBLP:conf/eurocrypt/Stadler96}
Markus Stadler.
\newblock Publicly verifiable secret sharing.
\newblock In Maurer \cite{DBLP:conf/eurocrypt/96}, pages 190--199.

\bibitem{gapdehn}
Chris Staecker.
\newblock dehn.gap.
\newblock http://cstaecker.fairfield.edu/~cstaecker/files/gap/dehn.gap.

\end{thebibliography}

\end{document}